\documentclass{article}
\usepackage[english]{babel}
\usepackage[T1]{fontenc}
\usepackage{amsmath, amsfonts, amsthm, amssymb, graphics}

\newcommand{\e}{\varepsilon}
\DeclareMathOperator{\p}{P}
\newcommand{\E}{\mathrm{E}}
\newcommand{\Var}{\mathrm{Var}}
\DeclareMathOperator{\BIC}{BIC}
\DeclareMathOperator{\mBIC}{mBIC}
\DeclareMathOperator{\mBIK}{mBIC2}
\DeclareMathOperator{\EBIC}{EBIC}
\DeclareMathOperator{\RSS}{RSS}
\DeclareMathOperator{\Tr}{Tr}
\newtheorem{theorem}{Theorem}
\newtheorem{lemma}{Lemma}
\newtheorem{fact}{Fact}
\newtheorem{definition}{Definition}

\begin{document}

\begin{center}{\Large\bf Consistency of modified versions of Bayesian Information Criterion in sparse linear regression with subgaussian errors}
\vspace{0.2in}

\noindent{\large Piotr Szulc}
\vspace{0.2in}

\noindent {\it \small University of Wroc\l{}aw, Poland}
\end{center}

\begin{abstract}
We consider a sparse linear regression model, when the number of available predictors, \(p\), is much larger than the sample size, \(n\), and the number of non-zero coefficients, \(p_0\), is small. To choose the regression model in this situation, we cannot use classical model selection criteria. In recent years, special methods have been proposed to deal with this type of problem, for example modified versions of Bayesian Information Criterion, like mBIC or mBIC2. It was shown that these criteria are consistent under the assumption that both \(n\) and \(p\) as well as \(p_0\) tend to infinity and the error term is normally distributed \cite{szulc}. In this article we prove the consistency of mBIC and mBIC2 under the assumption that the error term is a subgaussian random variable.
\end{abstract}

\section{Introduction}
One of purposes of analysis of large data sets is to determine which explanatory variables have a significant relationship with an explained variable. To achieve this goal in the context of the linear regression, it is natural to consider classical criteria like Akaike Information Criterion (AIC, Akaike 1974) or Bayesian Information Criterion (BIC, Schwarz 1978). However, it is now well known that this is not a good idea when the total number of available predictors, \(p\), is comparable or larger than a number of observations, \(n\). Specifically, Bogdan et al. (2008) showed that if \(\frac{p}{\sqrt{n\ln n}}\rightarrow \infty\), then the expected number of false predictors detected by BIC may go to infinity. Since AIC for \(n\geq 8\) selects more regressors than BIC, false discove\-ries appear even more often.

Therefore, we need other criteria that would find correct models in the \textit{large \(p\), small \(n\)} problem. The construction of such criteria is possible in a situation where the data-generating model is \textit{sparse}, that is the number of true predictors, \(p_0\), is small. It turns out that in such a case we can modify BIC to suit our needs.

In this paper we prove a desirable property (the consistency) of two modifications of BIC, mBIC and mBIC2, in the situation when the error term in the linear regression model is gaussian or subgaussian. Luo and Chen \cite{ebic} proved that a simi\-lar modification of BIC, EBIC, is consistent, but that result cannot be transferred directly to mBIC and mBIC2. Besides, the proof only applies to the situation when the error term is gaussian.

We would like to mention that one can find different approaches to a selection of variables in the \textit{large \(p\), small \(n\)} problem. For example, we can first perform so called \textit{screening}, in which we do some simple tests in order to remove most predictors, so that we end up with a traditional setting when \(n>p\). Pokarowski et al. (2015) showed that when we limit a number of predictors in this way (using LASSO) and then apply so called Generalized Information Criterion, the whole procedure is consistent.  

\section{mBIC and mBIC2}
Let \(\textbf{Y} = (Y_1,\ldots,Y_n)\) be a vector of random variables connected with values of the deterministic matrix,
\begin{equation*}
\textbf{X} = 
 \begin{pmatrix}
  1 & x_{1,1} & \cdots & x_{1,p} \\
  1 & x_{2,1} & \cdots & x_{2,p} \\
  \vdots & \vdots & \ddots & \vdots \\
  1 & x_{n,1} & \cdots & x_{n,p} 
 \end{pmatrix},
\end{equation*}
in the following way:
\begin{equation}
\label{reg}
Y_i = \beta_0 + \sum_{j=1}^p \beta_j x_{ij} + \e_i.
\end{equation}
where \(\boldsymbol\e = (\e_1,\ldots,\e_n)\) is a vector of independent and identically distributed random variables with mean zero and variance \(\sigma^2\) and \(\boldsymbol\beta=(\beta_1,\ldots,\beta_p)\) are unknown parameters.

We denote by \(s\) any subset of \(\{1,\ldots,p\}\) and by \(k(s)\) its size. With this notation, the expression that we minimize in BIC can be written as
\begin{equation}
\label{bic}
\BIC(s)= n\ln\RSS(s)+k(s)\ln n,
\end{equation}
where \(\RSS(s)\) is the residual sum of squares for a model \(s\). BIC was derived in the Bayesian context and assumes equal prior probability for all models. In the results, the prior distribution on the size of the true model is binomial, \(\mathcal{B}(p,1/2)\). Because this distribution concentrates almost entirely on \([p/2-3\sqrt{p},~p/2+3\sqrt{p}]\), it does not agree with the sparsity assumption.

The above observation leads to a natural modification of BIC by replacing the uniform distribution with a different one, more in line with the expectations that the number of true predictors is small. In 2004 Bogdan et al. \cite{mbic} proposed a modification called mBIC, in which the prior distribution on \(p_0\) is \(\mathcal{B}(p,c/p)\), where \(c\) is the expected value of the size of the true model. The resulting formula is
\begin{equation*}
\mBIC(s) = n\ln\RSS(s)+k(s)\ln n+2k(s)\ln\left(\frac{p}{c}-1\right).
\end{equation*}
Note that taking \(c=p/2\) the last component disappears and mBIC reduces to BIC. If we do not have any expectations about true size, it was shown that \(c=4\) is a good choice, that is the overall
type I error (FWER) is controlled at the level below 10\% if \(n\geq 200\), \(p\geq 20\) and \(\textbf{X}\) is orthogonal. 

If FWER is not our priority and we prefer to focus on the false discovery rate (FDR), mBIC2 \cite{mbic2} is a better choice. In this criterion we minimize
\begin{equation*}
\mBIK(s) = n\ln \RSS(s) + k(s)\ln n + 2k(s)\ln\left(\frac{p}{c}-1\right) -2\ln (k(s)!).
\end{equation*}
As shown in \cite{mbic2}, if we choose \(c=4\), FDR is controlled at the level below 10\% if \(n\geq 200\), \(p\geq 20\) and \(\textbf{X}\) is orthogonal.

The idea of mBIC was extended by Chen \cite{ebic2} in constructing EBIC, which in its standard version uses the uniform prior on \(p_0\). More generally, the expression minimized by EBIC can be written as
\begin{equation*}
\EBIC(s) = n\ln\RSS(s) + k(s)\ln n + 2\gamma\ln \binom{p}{k(s)},
\end{equation*}  
where \(\gamma\geq 0\). Luo and Chen \cite{ebic} proved that EBIC is consistent when \(n\), \(p\) and \(p_0\) go to infinity. It was shown in \cite{zak} that, in case of sparse models, EBIC and mBIC2 are asymptotically equivalent. However, the prior distribution in EBIC does not match the sparsity assumption very well because still the expected value of the size of the true model is \(p/2\).

\section{Consistency}
Consider the model (\ref{reg}) and let \(s_0\) be the true set of predictors, that is \(s_0=\{j\colon~\beta_j\neq 0,~j\in \{1,\ldots,p\}\}\). We say that a criterion \(C\) is consistent when
\begin{equation*}
\p\left(\underset{\substack{s\colon k(s)\leq K_n \\ s \neq s_0}}{\forall} C(s)>C(s_0)\right)\xrightarrow{n \rightarrow \infty} 1,
\end{equation*}
where \(K_n\) is the maximum size of searched models. This size for a fixed \(n\) has to be limited because when \(p\geq n\), it usually occurs \(\min\{\RSS(s)\colon k(s)=n\}=0\). However, in our consideration we allow \(K_n\) to go to infinity, thanks to which we are able to identify the correct model of any size for a sufficiently large \(n\).

We present proofs of the conistency of mBIC and mBIC2 in two versions: in case when the errors \(\e_i\) are gaussian and in more general case (with strongest assumptions), when they are subgaussian. The first situation was considered in \cite{szulc} but now the assumptions are weakened.

\subsection{Notations and basic assumptions}
Denote by \(X(s)\) a matrix composed of columns of \(X\) with indexes in \(s\) and, analogically, \(\beta(s)\) is a vector with elements of \(\beta\) with indexes in \(s\). Let \(H_n(s)\) be a matrix of the orthogonal projection on the space spanned by columns of \(X(s)\), that is \(H(s)=X(s)[X(s)^TX(s)]^{-1}X(s)^T\). Next, denote \(\Delta(s) = ||\mu-H(s)\mu||^2\), where \(\mu=X(s_0)\beta(s_0)\). Recall that we have \(p+1\) columns in \(X\) and the size of the true model is \(p_0\). It should be remembered that the parameters given above (in particular \(p\) and \(p_0\)) depend on \(n\). To simplify notation, we do not signal this as well in the case of \(K_n\), writing \(K\). We should also emphasize that because we are considering the model (\ref{reg}), we assume that the matrix \(X\) is deterministic.

We consider situation \(p\geq n\), in which some columns in \(X\) can be represented
as a linear combination of others. Hence, a vector of the expected values of \(Y\) can usually be represented by many combinations of available predictors, and as a result, one cannot talk about a single correct model. Therefore, we need a condition guaranteeing identification of the true model. Articles \cite{ident} and \cite{ident2} present appropriate assumptions in terms of the matrix \(X\) and the magnitude of the actual regression coefficients. Here we use the weaker and more convenient condition from \cite{ebic}, which is expressed in the language of \(\Delta(s)\).

\textbf{Consistency condition}
\begin{equation}
\label{Ident}
\lim_{n\rightarrow \infty}\min\left\{\frac{\Delta(s)}{p_0\ln p}\colon s_0 \nsubseteq s, k(s)\leq K\right\}=\infty,
\end{equation}
where \(K=kp_0\) for any fixed \(k>1\).

Such a condition implies that when \(X^TX=nI\), where \(I\) is the identity matrix of size \(n\),
\begin{equation}
\label{bety}
\sqrt{\frac{n}{p_0\ln p}}\min_{j\in s_0}|\beta_j| \rightarrow\infty.
\end{equation}
To show that, denote by \(l\) an element of \(s_0\) and by \(s_{-l}\) the set \(s_0\) without \(l\). We have
\begin{equation*}
\Delta(s_{-l}) = ||[I-H(s_{-l})]X(s_0)\beta(s_0)||^2 = ||[I-H(s_{-l})]X(\{l\})\beta_l||^2 = n\beta_l^2.
\end{equation*}
Hence,
\begin{equation*}
\min\left\{\frac{\Delta(s)}{p_0\ln p}\colon s_0 \nsubseteq s, k(s)\leq K\right\} \leq \frac{\Delta(s_{-l})}{p_0\ln p} \leq \frac{n}{p_0\ln p}\beta_l^2.
\end{equation*}
Because the above inequality holds for any \(l\in s_0\), the condition (\ref{Ident}) implies (\ref{bety}). So we allow the coefficients from the true model to go to zero, but not arbitrarily fast. Furthermore, as shown in \cite{ebic}, (\ref{bety}) implies (\ref{Ident}) if additionally the sparse Riesz condition holds:
\begin{align}
\label{riesz}
0 < c_1 & \leq \min_{s\colon k(s)\leq K}\left\{\lambda_1\left(\frac{1}{n}X(s)^T X(s)\right)\right\}\nonumber\\
& \leq \max_{s\colon k(s)\leq K}\left\{\lambda_2\left(\frac{1}{n}X(s)^T X(s)\right)\right\}\leq c_2 < \infty,
\end{align}
where \(\lambda_1\) and \(\lambda_2\) are the smallest and the largest eigenvalues, respectively.

\subsection{Gaussian error}
We consider the following criteria:
\begin{align*}
\mBIC(s) & = n\ln\RSS(s) + k(s)\ln n + 2k(s)\ln p,\\
\mBIK(s) & = n\ln\RSS(s_0) + k(s)\ln n + 2k(s)\ln p -2\ln(k(s)!).
\end{align*}
To simplify calculations, we replaced the expression \(\ln(p/c-1)\) by \(\ln p\), which does not affect the asymptotic properties of the criteria.

We begin with proving the auxiliary lemma associated with the \(\chi^2\) distribution, next we formulate and prove theorems about the consistency of mBIC and mBIC2. In both cases we assume that \(p_0\ln p=O(n^{\alpha})\), where \(\alpha\in [\frac{1}{2},1)\).  The higher \(\alpha\), the stronger the additional assumption on the maximum \(p\) is needed. We will show later that just after a slight modification of mBIC and mBIC2, it is enough to assume \(p_0\ln p=o(n)\). 

\begin{lemma}
\label{lemat_add}
Let \(\chi^2_j\) be a random variable with \(\chi^2\) distribution with \(j\) degrees of freedom. Denote \(m_j = 2j\left(\ln p + \frac{1}{2}\ln\ln p\right)\). If \(p\rightarrow \infty\), then
\begin{equation*}
\sum_{j=1}^K \binom{p}{j}\p(\chi^2_j\geq m_j) \rightarrow 0,
\end{equation*}
where \(K\leq p\) (but it can also go to infinity).
\end{lemma}

\begin{proof}
Using Lemma 2 from \cite{mielniczuk}, we can write
\begin{equation*}
\p(\chi^2_j\geq m) < 2\frac{(m/2)^{j/2-1}e^{-m/2}}{\Gamma(j/2)}.
\end{equation*}
We have \(\binom{p}{j} \leq \frac{p^j}{j!}\). Because \(\frac{j!\Gamma(j/2)}{(2j)^{j/2}} \rightarrow\infty\),
there is \(c>0\) that \(cj!\Gamma(j/2) > (2j)^{j/2}\). Therefore,
\begin{align*}
\binom{p}{j} \p(\chi_j^2\geq m_j) & \leq 2 \frac{p^j}{j!} \frac{\left(j\ln p + \frac{1}{2}j\ln\ln p\right)^{j/2-1}} {\Gamma(j/2)\exp\left(j\ln p + \frac{1}{2}j\ln\ln p\right)}\\
& < \frac{2c}{\ln p} \frac{j^{j/2}\left(\ln p + \frac{1}{2}\ln\ln p\right)^{j/2}}{(2j)^{j/2}(\ln p)^{j/2}}\\
& = \frac{2c}{\ln p} \left[\sqrt{\frac{\ln p + \frac{1}{2}\ln\ln p}{2\ln p}}\right]^j.
\end{align*}

The base of the power in the last expression can be limited by a constant \(q<1\), so we can write
\begin{align*}
\sum_{j=1}^K\binom{p}{j}\p(\chi^2_j\geq m_j) & < \frac{2c}{\ln p} \sum_{j=1}^K q^j < \frac{2c}{\ln p} \frac{q}{1-q} \xrightarrow{p\rightarrow \infty} 0,
\end{align*}
which proves the lemma.
\end{proof}

Before we go to the main theorems, we will prove a simple fact about \(\chi^2\) distribution. First, let us introduce additional notation. By \(o_p(a_n)\), where \(a_n\) is a numerical sequence, we understand a random variable \(X_n\), for which the quotient \(X_n/a_n\) converges in probability to zero when \(n\) goes to infinity, \textit{i.e.} for any \(\delta > 0\), \(\lim_{n\rightarrow\infty} \p\left(\left|X_n/a_n\right| > \delta\right) = 0\). Furthermore, we write \(X_n = O_p(a_n)\) if for any \(\delta>0\) exists \(M>0\), that for every \(n\) we have \(\p\left(\left|X_n/a_n\right| > M\right) < \delta\).

\begin{fact}\label{lem3}
Let \(\chi^2_j\) be a random variable with \(\chi^2\) distribution and \(j\) degrees of freedom. Then \(\chi^2_j = j(1+o_p(1))\). Besides, \(o_p(1)=O_p\left(1/\sqrt j\right)\).
\end{fact}
\begin{proof}
We have \(\E(\chi^2_j) = j\) and \(\Var(\chi^2_j) = 2j\). Using Chebyshev's inequality, for any \(\delta\) the following inequality holds:
\begin{equation*}
\p\left(|\chi^2_j-j|\geq \delta\sqrt{2j}\right) \leq \frac{1}{\delta^2}.
\end{equation*}
Hence \(\chi^2_j = j + O_p(\sqrt{j})\), which proves the fact.
\end{proof}

\begin{theorem}
Assume the model (\ref{reg}), where \(\e_i\) are gaussian, \(p\geq n\) and let the condition \ref{Ident} hold. If \(p_0\ln p=O(n^{\alpha})\) and \(\ln p=o(n^{1-\alpha}\ln n)\), where \(\alpha\in \left[\frac{1}{2},1\right)\), then mBIC is consistent, that is
\begin{equation*}
\p\left(\underset{\substack{s\colon k(s)\leq K \\ s \neq s_0}}{\forall} \mBIC(s)>\mBIC(s_0)\right)\xrightarrow{n \rightarrow \infty} 1,
\end{equation*}
where \(K=kp_0\) for any fixed \(k>1\).
\end{theorem}

\begin{proof}
Denote by \(s_0\) the true model and by \(s\neq s_0\) any different model. We will show that when \(n\) is large enough, the difference \(\mBIC(s)-\mBIC(s_0)\) is larger than zero with probability going to one.

Note that if \(\e_i\) is \(\mathcal{N}(0,\sigma^2)\), we can present it as \(\e_i=\sigma\e_i'\), where \(\e_i'\) is \(\mathcal{N}(0,1)\). Therefore, the model \(y=X\beta+\e\) is equivalent to \(y'=X'\beta+\e'\), where \(y'=y/\sigma\) and \(X'=X/\sigma\). We have
\begin{equation}
\label{skal}
\frac{\RSS(s)}{\RSS(s_0)} = \frac{y^T[I-H(s)]y}{y^T[I-H(s_0)]y} = \frac{\sigma^2y'^T[I-H'(s)]y'}{\sigma^2y'^T[I-H'(s_0)]y'} = \frac{\RSS'(s)}{\RSS'(s_0)},
\end{equation}
so the quotient \(\RSS(s)/\RSS(s_0)\) does not depend on a scale of \(y\) and \(X\). Because the proof is based on an estimation of \(\RSS(s)/\RSS(s_0)\), we can assume \(\sigma^2=1\) without loss of generality.

We can write
\begin{align}
\label{roznica}
& \mBIC(s)-\mBIC(s_0) \nonumber\\
& = n\ln\frac{\RSS(s)}{\RSS(s_0)} + (k(s)-p_0)\ln n + 2(k(s)-p_0)\ln p\nonumber\\ 
& \geq n\ln\frac{\RSS(s)}{\RSS(s_0)} - 3p_0\ln p =  n\ln\left(1+\frac{\RSS(s)-\RSS(s_0)}{\RSS(s_0)}\right)  - 3p_0\ln p.
\end{align}
Let us begin with \(\RSS(s_0)\). Note that the residual sum of squares can be written as \(\RSS(s) = y^T[I-H(s)]y\), where \(I\) is the identity matrix of the size \(n\). Because \([I-H(s_0)]X(s_0)=0\), we get
\begin{align*}
y^T[I-H(s_0)]y &= (X(s_0)\beta(s_0)+\e)^T[I-H(s_0)](X(s_0)\beta(s_0)+\e) =\\
&= \e^T[I-H(s_0)]\e.
\end{align*}
Because \(I-H(s_0)\) is an symmetric idempotent matrix with rank \(n-p_0\), \(\e^T[I-H(s_0)]\e\) has the \(\chi^2\) distribution with \(n-p_0\) degrees of freedom. Using the fact \ref{lem3}, we can write
\begin{equation*}
\e^T[I-H(s_0)]\e =(n-p_0)(1+o_p(1))= n(1+o_p(1)).
\end{equation*}
The last equality comes from the assumption that \(p_0=o(n)\).

Now let us estimate \(\RSS(s)-\RSS(s_0)\). First, assume that \(s\) does not include the true model, that is \(s_0 \not\subset s\). We have
\begin{align}
\label{suma}
\RSS(s)-\RSS(s_0) & = y^T[I-H(s)]y - \e^T[I-H(s_0)]\e \nonumber\\
& = \Delta(s)+2\mu^T[I-H(s)]\e + \e^T H(s_0)\e-\e^T H(s)\e,
\end{align}
what we get by substituting \(\mu+\e\) in place of \(y\) and using the inequality \(\mu^T[I-H(s)]\e = \e^T[I-H(s)]\mu\).

We will estimate components of the above sum. Using again the fact \ref{lem3}, we get
\begin{equation}
\label{suma1}
\e^TH(s_0)\e = p_0(1+o_p(1)).
\end{equation}

Let \(m_K = 2K\left(\ln p + \frac{1}{2}\ln\ln p\right)\). Because the random variable \(\e^TH(s)\e\) is \(\chi^2\), we will denote it as \(\chi_j^2\), where \(j=k(s)\). From the Bonferroni inequality, we get
\begin{align*}
\p\left(\underset{\substack{s\colon k(s)\leq K \\ s_0 \not\subset s}}{\exists}\e^T H(s)\e \geq m_K\right) = \p\left(\underset{1\leq j \leq K}{\exists}~\underset{\substack{s\colon k(s)=j \\ s_0 \not\subset s}}{\exists}\e^T H(s)\e \geq m_K\right)\\
\leq \sum_{j=1}^K\binom{p}{j}\p(\chi_j^2\geq m_K) \leq \sum_{j=1}^K\binom{p}{j}\p(\chi_j^2\geq m_j).
\end{align*}
From the lemma \ref{lemat_add} the last sum goes to zero, so
\begin{equation}
\label{suma2}
\underset{\substack{s\colon k(s)\leq K \\ s_0 \not\subset s}}{\forall}~\e^TH(s)\e \leq O_p(K \ln{p}),
\end{equation}
where \(O_p(K \ln{p})\) does not depend on \(s\).

Now, we show that
\begin{equation}
\label{suma3}
\underset{\substack{s\colon k(s)\leq K \\ s_0 \not\subset s}}{\forall}~|\mu^T[I-H(s)]\e| = \sqrt{\Delta(s)O_p(K\ln{p})},
\end{equation}
where \(O_p(K \ln{p})\) does not depend on \(s\). We can write
\begin{equation*}
\mu^T[I-H(s)]\e=\sqrt{\Delta(s)}\e'(s),
\end{equation*}
where \(\e'(s)\sim \mathcal{N}(0,1)\), because
\begin{align*}
\E(\mu^T[I-H(s)]\e) & = \mu^T[I-H(s)]\E(\e)=0=\E\left(\sqrt{\Delta(s)}\e'(s)\right)\\
\Var(\mu^T[I-H(s)]\e) & = \mu^T[I-H(s)]\Var(\e)(\mu^T[I-H(s)])^T\\
& = \mu^T[I-H(s)]\mu = \Var\left(\sqrt{\Delta(s)}\e'(s)\right).
\end{align*}
For \(m_K\) defined above we have
\begin{align*}
\p\left(\underset{\substack{s\colon k(s)\leq K \\ s_0 \not\subset s}}{\exists}\,|\e'(s)| \geq \sqrt m_K\right) \leq \sum_{j=1}^K \binom{p}{j} \p\left(|\e'(s)|\geq \sqrt m_K\right)\\
= \sum_{j=1}^K \binom{p}{j}\p(\chi_1^2(s)\geq m_K) \leq \sum_{j=1}^K \binom{p}{j}\p(\chi_j^2(s)\geq m_K),
\end{align*}
because \(\p(\chi_1^2\geq m_K)\leq \p(\chi_j^2\geq m_K)\). Using again the lemma \ref{lemat_add}, we get
\begin{equation*}
\underset{\substack{s\colon k(s)\leq K \\ s_0 \not\subset s}}{\forall}~|\mu^T[I-H(s)]\e| = \sqrt{\Delta(s)O_p(K\ln{p})}.
\end{equation*}

Let us go back to (\ref{suma}). Using (\ref{suma1}), (\ref{suma2}), (\ref{suma3}) and the condidtion( \ref{Ident}), we can write
\begin{align}
\label{RSSroznica}
& \underset{\substack{s\colon k(s)\leq K \\ s_0 \not\subset s}}{\forall}~\RSS(s)-\RSS(s_0) \nonumber\\
& ~~~~ = \Delta(s)\left(1 + \frac{2\mu^T[I-H(s)]\e}{\Delta(s)} + \frac{\e^T H(s_0)\e}{\Delta(s)} - \frac{\e^T H(s)\e}{\Delta(s)}\right) \nonumber\\
& ~~~~ = \Delta(s)(1+o_p(1)),
\end{align}
where \(o_p(1)\) converges in probability to zero and does not depend on \(s\). Now we can go further in (\ref{roznica}): for any constant \(C>0\) and \(n\) large enough, we have
\begin{align}
\label{norm_nier}
& \mBIC(s)-\mBIC(s_0) =\nonumber\\
&~~ =n\ln\left(1+\frac{\Delta(s)}{n}(1+o_p(1))\right) - 3p_0\ln p\nonumber\\
&~~ \geq n\ln\left(1+\frac{Cp_0\ln p}{n}(1+o_p(1))\right)-3p_0\ln p\nonumber\\
&~~ = Cp_0\ln p\ln\left(1+\frac{Cp_0\ln p}{n}(1+o_p(1))\right)^{\frac{n}{Cp_0\ln p}}-3p_0\ln p,
\end{align}
The expression
\begin{equation*}
\left(1+\frac{Cp_0\ln p}{n}(1+o_p(1))\right)^{\frac{n}{Cp_0\ln p}}
\end{equation*}
goes to \(e\), because \(p_0\ln p=o(n)\) and for every \(k\rightarrow\infty\) we have
\begin{equation*}
\left(1+\frac{1}{k}(1+o_p(1))\right)^k \xrightarrow{P} e.
\end{equation*}
Therefore, if \(C\) is large enough, the above difference is larger than zero for every \(s\) such that \(k(s)\leq K\) and \(s_0 \not\subset s\) with probability going to one.

\begin{center} * * * \end{center}

Consider the case \(s_0\subset s\). We have \([I-H(s)]X(s_0)=0\), so
\begin{equation*}
\RSS(s) =  y^T[I-H(s)]y = \e^T[I-H(s)]\e
\end{equation*}
and
\begin{align*}
\RSS(s_0)-\RSS(s) & = \e^T[I-H(s_0)]\e - \e^T[I-H(s)]\e \\
& = \e^T[H(s)-H(s_0)]\e = \chi_j^2(s),
\end{align*}
where \(j=k(s)-p_0\). In that case, we can write
\begin{align*}
\ln\frac{\RSS(s_0)}{\RSS(s)} &=\ln\left(1+\frac{\RSS(s_0)-\RSS(s)}{\RSS(s_0)-[\RSS(s_0)-\RSS(s)]}\right)\\
& \leq \frac{\chi_j^2(s)}{\e^T[I-H(s_0)]\e-\chi_j^2(s)},
\end{align*}
what results from \(\ln(1+x)\leq x\).

Using the fact \ref{lem3}, we get
\begin{equation*}
\e^T[I-H(s_0)]\e=n + O_p(\sqrt n).
\end{equation*}

We have
\begin{align}
\label{subt}
& \p\left(\underset{\substack{s\colon k(s)\leq K \\ s_0 \subset s,\,s\neq s_0}}{\forall}~\mBIC(s)-\mBIC(s_0) > 0\right)\nonumber\\
&~~~~= \p\left(\underset{\substack{s\colon k(s)\leq K \\ s_0 \subset s,\,s\neq s_0}}{\forall}~ \ln\frac{\RSS(s_0)}{\RSS(s)} < \frac{j\ln n + 2j\ln p}{n}\right)\nonumber\\
&~~~~\geq \p\left(\underset{\substack{s\colon k(s)\leq K \\ s_0 \subset s,\,s\neq s_0}}{\forall}~\frac{\chi_j^2(s)}{n + O_p(\sqrt n)-\chi_j^2(s)} < \frac{j\ln n +2j\ln p}{n}\right)\nonumber\\
&~~~~= 1-\p\left(\underset{\substack{s\colon k(s)\leq K \\ s_0 \subset s,\,s\neq s_0}}{\exists}~\chi_j^2(s) \geq \frac{n+O_p(\sqrt n)}{n + j\ln n + 2j\ln p}(j\ln n +2j\ln p)\right)\nonumber\\
&~~~~= 1-\p\left(\underset{\substack{s\colon k(s)\leq K \\ s_0 \subset s,\,s\neq s_0}}{\exists}~\chi_j^2(s) \geq (1+O_p(n^{\alpha-1}))(j\ln n + 2j\ln p)\right),
\end{align}
where \(O_p(n^{\alpha-1})\) does not depend on \(s\). The last equality comes from the assumption that \(p_0\ln p = O(n^{\alpha})\), because then, if \(n\) is large enough, we have
\begin{equation*}
\frac{n+O_p(\sqrt n)}{n + j\ln n + 2j\ln p} = \frac{n+O_p(\sqrt n)}{n + O(n^{\alpha})} = 1+\frac{O_p(\sqrt n)-O(n^{\alpha})}{n+O(n^{\alpha})} = 1+O_p(n^{\alpha-1}).
\end{equation*}

Let \(m_j = 2j\left(\ln p + \frac{1}{2}\ln\ln p\right)\). From the Bonferroni inequality and the lemma \ref{lemat_add} we get
\begin{align}
\label{mj}
\p\left(\underset{1\leq j \leq K-p_0}{\exists}~\underset{\substack{s\colon k(s)=j+p_0 \\ s_0 \subset s}}{\exists}\,\chi^2_j(s) \geq m_j\right) & \leq \sum_{j=1}^{K-p_0}\binom{p-p_0}{j} \p(\chi_j^2(s)\geq m_j)\nonumber\\
& \leq \sum_{j=1}^K \binom{p}{j} \p(\chi_j^2(s)\geq m_j) \xrightarrow{n\rightarrow \infty} 0.
\end{align} 

We will show that with probability going to one we have
\begin{equation}
\label{nier1}
(1+O_p(n^{\alpha-1}))(j\ln n + 2j\ln p) > m_j.
\end{equation}
Note that the left side of the inequality can be written as
\begin{equation*}
j\ln n + 2j\ln p + O_p(n^{\alpha-1})j\ln p.
\end{equation*}
From the assumption \(\ln p = o(n^{1-\alpha}\ln n)\) it follows that
\begin{equation*}
\ln n > \ln \ln p + O_p(n^{\alpha-1})\ln p
\end{equation*}
for sufficiently large \(n\), and this implies (\ref{nier1}). In that case, from (\ref{subt}) and (\ref{mj}) we get
\begin{equation*}
\p\left(\underset{\substack{s\colon k(s)\leq K \\ s_0 \subset s,\,s\neq s_0}}{\forall}~\mBIC(s)-\mBIC(s_0) > 0\right) \xrightarrow{n \rightarrow \infty} 1. \qedhere
\end{equation*}
\end{proof}

\begin{center} * * * \end{center}

Now we formulate the analogous theorem about mBIC2.

\begin{theorem}
Assume the model \ref{reg}, where \(\e_i\) are gaussian, \(p\geq n\) and let the condition (\ref{Ident}) hold. If \(p_0^2\ln p=O(n^{\alpha})\) and \(\ln p=o(n^{1-\alpha}\ln n)\), where \(\alpha\in \left[\frac{1}{2},1\right)\), then mBIC2 is consistent, that is
\begin{equation*}
\p\left(\underset{\substack{s\colon k(s)\leq K \\ s \neq s_0}}{\forall} \mBIK(s)>\mBIK(s_0)\right)\xrightarrow{n \rightarrow \infty} 1,
\end{equation*}
where \(K=kp_0\) for any fixed \(k>1\).
\end{theorem}

\begin{proof}
Assume \(s_0 \not\subset s\). Using estimation from the mBIC proof, we can write
\begin{align}
\label{subg3}
& \mBIK(s)-\mBIK(s_0)\nonumber\\
& = \mBIC(s) - \mBIC(s_0) + 2\ln (p_0!) - 2\ln (k(s)!)\nonumber\\
& \geq Cp_0\ln p \ln\left(1 + \frac{Cp_0\ln p}{n}(1+o_p(1))\right)^{\frac{n}{Cp_0\ln p}} -(2k+3)p_0\ln p, 
\end{align}
where \(o_p(1)\) does not depend on \(s\), because
\begin{equation*}
\ln (k(s)!)\leq \ln ((kp_0)!)\leq kp_0\ln kp_0 \leq kp_0\ln p.
\end{equation*}
Analogically to the mBIC proof, for \(n\) large enough the difference \(\mBIK(s)-\mBIK(s_0)\) is larger than zero with the probability going to one.

Now, let \(s_0 \subset s\). We have to estimate \(\ln(k(s)!)-\ln(p_0!)\) more carefully. Because \(\ln(n!) \geq n\ln n -n\) and \(\ln(n!) \leq n\ln n\), so remembering that \(j=k(s)-p_0\), we can write
\begin{align*}
\ln (k(s)!) -\ln (p_0!) & = \ln((j+p_0)!) -\ln (p_0!) = \ln(1+p_0) +\ldots + \ln(j+p_0)\\
& \leq j\ln(j+p_0) \leq j\ln(kp_0).
\end{align*}

Therefore, we have
\begin{align}
\label{subt2}
& \p\left(\underset{\substack{s\colon k(s)\leq K \\ s_0 \subset s,\,s\neq s_0}}{\forall}~\mBIK(s)-\mBIK(s_0) > 0\right)\nonumber\\
& = \p\left(\underset{\substack{s\colon k(s)\leq K \\ s_0 \subset s,\,s\neq s_0}}{\forall}~\ln\frac{\RSS(s_0)}{\RSS(s)} < \frac{j\ln n +2j\ln p - 2j\ln(kp_0)}{n}\right)\nonumber\\
& = 1-\p\left(\underset{\substack{s\colon k(s)\leq K \\ s_0 \subset s,\,s\neq s_0}}{\exists}~\chi_j^2(s) \geq (1+O_p(n^{\alpha-1}))(j\ln n + 2j\ln p - 2j\ln(kp_0))\right),
\end{align}
where \(O_p(n^{\alpha-1})\) does not depend on \(s\).

Let \(m_j = 2j\left(\ln p + \frac{1}{2}\ln\ln p\right)\). We will show that with probability going to one
\begin{equation}
\label{nier2}
(1+O_p(n^{\alpha-1}))(j\ln n + 2j\ln p - 2j\ln(kp_0)) > m_j.
\end{equation}
Note that the left side of the inequality can be written as
\begin{equation*}
j\ln n + 2j\ln p - 2j\ln(kp_0) + O_p(n^{\alpha-1})j\ln p.
\end{equation*}
From the assumption \(\ln p = o(n^{1-\alpha}\ln n)\) i \(p_0^2\ln p = O(n^{\alpha})\), it follows that 
\begin{equation*}
\ln n > \ln(k^2p_0^2\ln p) + O_p(n^{\alpha-1})\ln p
\end{equation*}
for sufficiently large \(n\), and that implies (\ref{nier2}). In that case, we get the thesis from (\ref{mj}).
\end{proof}

\subsection{Extensions of mBIC and mBIC2}
The assumptions of the above theorems can be weakened if we consider the criteria in the following form:
\begin{align*}
\mBIC_{\gamma}(s) & = n\ln\RSS(s) + k(s)\ln n + 2\gamma k(s)\ln p,\\
\mBIK_{\gamma}(s) & = n\ln\RSS(s_0) + k(s)\ln n + 2\gamma k(s)\ln p -2\ln(k(s)!),
\end{align*}
where \(\gamma\geq 1\) is a constant.

\begin{theorem}
Assume the model (\ref{reg}), where \(\e_i\) are gaussian, \(p\geq n\) and let the condition (\ref{Ident}) hold. If \(p_0\ln p =o(n)\), than for \(\gamma>1\) \(\mBIC_{\gamma}\) is consistent.
\end{theorem}
\begin{proof}
The proof is analogous to mBIC. The only place where the assumption \(p_0\ln p =o(n)\) is not sufficient is (\ref{subt}). We have
\begin{align}
\label{subtgamma}
& \p\left(\underset{\substack{s\colon k(s)\leq K \\ s_0 \subset s,\,s\neq s_0}}{\forall}~\mBIC_{\gamma}(s)-\mBIC_{\gamma}(s_0) > 0\right)\nonumber\\
&~~~~= 1-\p\left(\underset{\substack{s\colon k(s)\leq K \\ s_0 \subset s,\,s\neq s_0}}{\exists}~\chi_j^2(s) \geq \frac{n+O_p(\sqrt n)}{n + j\ln n + 2\gamma j\ln p}(j\ln n +2\gamma j\ln p)\right)\nonumber\\
&~~~~= 1-\p\left(\underset{\substack{s\colon k(s)\leq K \\ s_0 \subset s,\,s\neq s_0}}{\exists}~\chi_j^2(s) \geq (1+o_p(1))(j\ln n + 2\gamma j\ln p)\right),
\end{align}
where \(o_p(1)\) does not depend on \(s\). The last equality comes from \(p_0\ln p =o(n)\). For \(m_j\) as in (\ref{mj}), with probability going to one we have
\begin{equation*}
(1+o_p(1))(j\ln n + 2\gamma j\ln p) > m_j,
\end{equation*}
if \(\gamma>1\), and then we get the thesis.
\end{proof}

\begin{theorem}
Assume the model (\ref{reg}), where \(\e_i\) are gaussian, \(p\geq n\) and let the condition (\ref{Ident}) hold. If \(p_0\ln p =o(n)\), than for \(\gamma>1\) \(\mBIK_{\gamma}\) is consistent.
\end{theorem}
\begin{proof}
The proof is analogous to mBIC2. The only place where the assumption \(p_0\ln p =o(n)\) is not sufficient is (\ref{subt2}). We have
\begin{align}
\label{subt4}
& \p\left(\underset{\substack{s\colon k(s)\leq K \\ s_0 \subset s,\,s\neq s_0}}{\forall}~\mBIK_{\gamma}(s)-\mBIK(s_0) > 0\right)\nonumber\\
& ~~~~ = 1-\p\left(\underset{\substack{s\colon k(s)\leq K \\ s_0 \subset s,\,s\neq s_0}}{\exists}~\chi_j^2(s) \geq (1+o_p(1))(j\ln n + 2\gamma j\ln p - 2j\ln(kp_0)\right),
\end{align}
where \(o_p(1)\) does not depend on \(s\). The last equality comes from \(p_0\ln p =o(n)\). For \(m_j\) as in (\ref{mj}), with probability going to one we have
\begin{equation*}
(1+o_p(1))(j\ln n + 2\gamma j\ln p - 2j\ln(kp_0)) > m_j,
\end{equation*}
if \(\gamma>1\), and then we have the thesis.
\end{proof}

\subsection{Subgaussian error}
mBIC and mBIC2 were constructed assuming that the random error is gaussian. In case of real data analysis, this assumption is often not met. We will show that when \(\e_i\) is subgaussian, extensions  \(\mBIC_{\gamma}\) and \(\mBIK_{\gamma}\) defined in the previous section are consistent. Note that in case of many distributions occurring in nature, we can limit the support (for example the growth of man can be neither negative nor greater than a certain number), and each distribution with the limited support is subgaussian.

The penalty \(2 \ln p \) in both criteria was chosen to eliminate false discoveries in the case of the normal distribution, but at the same time to retain the highest possible power. If we compare a gaussian variable with a subgaussian one with the same variance \(\sigma^2\), the latter can have much heavier tails. This fact suggests that the penalty \(2\ln p\) may be insufficient and if we want to maintain the same convergence conditions, we need \(2\gamma\ln p\) for \(\gamma>1\). In the following theorems, we assume that \(\gamma\) is at least equal to the ratio \(b^2/\sigma^2\), where \(b\) is a subgaussian parameter. When \(b^2/\sigma^2 > 1\), original mBIC and mBIC2 are still consistent (which we will also show), but with much stronger restrictions on \(p\).

We will start by presenting a few basic facts related to the subgaussian distribution and we will prove two lemmas.

\begin{definition}
We say that \(\e\) is \(b\)-subgaussian if there is a positive constant \(b\) that for every \(t \in \mathbb{R}\) we have \(\E\left(e^{t\e_i}\right)\leq e^{b^2t^2/2}\).
\end{definition}
\begin{fact}
\label{fakt1}
If \(\e\) is \(b\)-subgaussian, then \(\E(\e)=0\) and \(\Var(\e)\leq b^2\).
\end{fact}
\begin{fact}
If \(\e\) is \(b\)-subgaussian, then for any \(\alpha\) a random variable \(\alpha\e\) is \(|\alpha|b\)-subgaussian.
\end{fact}
\begin{fact}
\label{fakt3}
If variables \(\e_i\) are independent and \(b_i\)-subgaussian, \(i=1,\ldots,n\), than \(\sum_{i=1}^n \e_i\) is \(\sqrt{\sum_{i=1}^n b_i^2}\)-subgaussian. 
\end{fact}
\begin{fact}
\label{fakt4}
If a random variable \(\e\) is subgaussian, then there is a positive constant \(c\) that for every \(m>0\) we have \(\p(|\e|\geq m) \leq 2e^{-cm^2}\).
\end{fact}

\begin{lemma}
\label{lemat2}
Let \(\e = (\e_1,\ldots,\e_n)\) be a vector of independent random variables \(b\)-subgaussian distribution, that is \(\E\left(e^{t\e_i}\right)\leq e^{b^2t^2/2}\) for every \(t\in\mathbb{R}\). Let \(A_j\) be an symmetric idempotent matrix with rank \(j\) and size \(n\). Denote \(m_j = 2b^2j\left(\ln p + \sqrt{2\ln p}\right)\). If \(p\rightarrow\infty\), then
\begin{equation*}
\sum_{j=1}^K \binom{p}{j}\p(\e^TA_j\e\geq m_j) \rightarrow 0,
\end{equation*}
where \(n\leq K\leq p\) (both \(n\) and \(K\) may go to infinity).
\end{lemma}

\begin{proof}
It was shown in \cite{ineq} that the following inequality regarding a quadratic form \(||A\e||^2\) holds for every \(t>0\):
\begin{equation}
\label{ineq1}
\p\left(||A\e||^2 > b^2\left(\Tr (A^TA) + 2\sqrt{\Tr ((A^TA)^2) t} + 2||A^TA||t\right)\right) \leq \exp(-t),
\end{equation}
where \(||A^TA||\) is the spectral norm of \(A^TA\). For a matrix \(A_j\) we have \(A_j^TA_j = A_j\), \(\Tr A_j = j\) and \(||A_j^TA_j||=1\) (the maximal eigenvalue \(A_j\)). Hence, the inequality (\ref{ineq1}) can be written in a simpler form:
\begin{equation*}
\p\left(\e^TA_j\e > b^2\left(j + 2\sqrt{jt} + 2t\right)\right) \leq \exp(-t).
\end{equation*}
Let \(m=b^2(j+2\sqrt{jt} + 2t)\), then
\begin{equation}
\label{ineq2}
\p\left(\e^TA_j\e > m\right) \leq \exp\left(-\frac{m}{2b^2} + \sqrt{\frac{mj}{2b^2}-\frac{j^2}{4}}\right) \leq \exp\left(-\frac{m}{2b^2} + \sqrt{\frac{mj}{2b^2}}\right).
\end{equation}

Using the inequality (\ref{ineq2}) for \(m_j=2b^2j\left(\ln p + \sqrt{2\ln p}\right)\), we get
\begin{align*}
& \sum_{j=1}^K \binom{p}{j} \p\left(\e^TA_j\e > m_j\right)\\
& ~~~~\leq \sum_{j=1}^K p^j \exp\left(-j\ln p-j\sqrt{2\ln p}+j\sqrt{\ln p +\sqrt{2\ln p}}\right)\\
& ~~~~= \sum_{j=1}^K\exp\left(-j\sqrt{2\ln p}+j\sqrt{\ln p +\sqrt{2\ln p}}\right)\leq \sum_{j=1}^K\exp\left(-\frac{j}{3}\sqrt{\ln p}\right).
\end{align*}
The last inequality is true if \(p\) is sufficiently large. We get a geometric series \(\sum_{j=1}^K\exp\left(-\frac{1}{3}\sqrt{\ln p}\right)^j\) that sums up to
\begin{equation*}
\exp\left(-\frac{1}{3}\sqrt{\ln p}\right)\frac{1-\exp\left(-\frac{K}{3}\sqrt{\ln p}\right)}{1-\exp\left(-\frac{1}{3}\sqrt{\ln p}\right)}
\end{equation*}
and converges to zero when \(p\) goes to infinity. 
\end{proof}

\begin{lemma}
\label{eAe}
Let \(A_j = [a_{ik}]_{n\times n}\) be an symmetric idempotent matrix with rank \(j\) and \(\e = (\e_1,\ldots,\e_n)\) be a vector of independent random variables with the following properties: \(\E(\e)=0\), \(\E(\e^2)=1\) and \(\E(\e^4)<\infty\). Then \(\e^TA_j\e = j(1+o_p(1))\). Besides, \(o_p(1)=O_p\left(1/\sqrt j\right)\).
\end{lemma}
\begin{proof}
To estimate \(\e^T A_j\e\), we use formulas on \(\E(\e^T A_j\e)\) and \(\Var(\e^T A_j\e)\) from \cite{Var} (the second one is true for \(\E\e=0\)).
\begin{align*}
\E(\e^T A_j\e) &= \Tr(A_j\Var(\e)) + (\E(\e))^TA_j\E(\e) = \Tr(A_j)=j,\\
\Var(\e^T A_j\e) &= [\E(\e^4)-3(\E(\e^2))^2]\sum_{i=1}^n a_{ii}^2 + [(\E(\e^2))^2-1](\Tr(A_j))^2 \\
& + 2(\E(\e^2))^2\Tr(A_j^2) \leq cj
\end{align*}
for a constant \(c\) because we assumed that \(\E(\e)=0\), \(\E(\e^2)=1\), \(\E(\e^4)<\infty\), \(\Tr(A_j^2)=\Tr(A_j)=j\) and
\begin{equation*}
\sum_{i=1}^n a_{ii}^2 \leq \sum_{i,k=1}^n a_{ik}^2 = \Tr (A_j^T A_j) = \Tr(A_j) =j.
\end{equation*}
Using Chebyshev's inequality, for any \(\delta\) we have
\begin{equation*}
\p\left(\left|\e^TA_j\e-j\right|\geq \delta\sqrt{cj}\right) \leq \p\left(\left|\e^TA_j\e-j\right|\geq \delta\sqrt{\Var(\e^T A_j\e)}\right) \leq \frac{1}{\delta^2}.
\end{equation*}
Hence \(\e^TA_j\e = j + O_p(\sqrt{j})\), which proves the lemma.
\end{proof}

\begin{theorem}
Assume the model (\ref{reg}), where \(\e_i\) are \textit{b}-subgaussian with the variance \(\sigma^2\), \(p\geq n\) and let the condition (\ref{Ident}) hold. If \(p_0\ln p=o(n)\), then for \(\gamma> b^2/\sigma^2\) \(\mBIC_{\gamma}\) is consistent, that
\begin{equation*}
\p\left(\underset{\substack{s\colon k(s)\leq K \\ s \neq s_0}}{\forall} \mBIC_{\gamma}(s)>\mBIC_{\gamma}(s_0)\right)\xrightarrow{n \rightarrow \infty} 1,
\end{equation*}
where \(K=kp_0\) for any fixed \(k>1\). If \(\gamma= b^2/\sigma^2\), \(\mBIC_{\gamma}\) is consistent as long as \(\sqrt{\ln p}=o(\ln n)\) and \(p_0=O(n^{\alpha})\) for \(\alpha<1\). If \(1\leq\gamma<b^2/\sigma^2\), besides \(p_0=O(n^{\alpha})\) it has to be \(p=O(n^{\delta})\), where \(\delta < \frac{1}{2(b^2/\sigma^2-\gamma)}\).
\end{theorem}

\begin{proof}
Note that if \(\e_i\) is \(b\)-subgaussian with the variance \(\sigma^2\), we can write \(\e_i=\sigma\e_i'\), where \(\e_i'\) is \(b/\sigma\)-subgaussian with the variance 1. Therefore, arguing as in \ref{skal}, we will consider \(\e_i\) with the variance 1 and subgaussian parameter \(B=b/\sigma\).

Let us begin with the case when a set \(s\) does not include the true model, that is \(s_0 \not\subset s\). We will estimate the difference \(\mBIC_{\gamma}(s)>\mBIC_{\gamma}(s_0)\):
\begin{align*}
& \mBIC_{\gamma}(s)-\mBIC_{\gamma}(s_0) \\
& ~~~~~~~~~~= n\ln\frac{\RSS(s)}{\RSS(s_0)} + (k(s)-p_0)\ln n + 2\gamma(k(s)-p_0)\ln p\\
& ~~~~~~~~~~\geq n\ln\left(1+\frac{\RSS(s)-\RSS(s_0)}{\RSS(s_0)}\right)  - (2\gamma+1)p_0\ln p.
\end{align*}
First, let us look at \(\RSS(s_0)\). We can write that \(\RSS(s_0) = y^T[I-H(s_0)]y = \e^T[I-H(s_0)]\e\). Because \(I-H(s_0)\) is a symmetric idempotent matrix with rank \(n-p_0\), so using the lemma \ref{eAe}, we get
\begin{equation}
\label{subg2}
\RSS(s_0) = (n-p_0)(1+o_p(1)) =n(1+o_p(1)).
\end{equation}

To estimate \(\RSS(s)-\RSS(s_0)\), we again present this difference as in (\ref{suma}):
\begin{equation}
\label{subg1}
\RSS(s)-\RSS(s_0) = \Delta(s) + 2\mu^T[I-H(s)]\e + \e^TH(s_0)\e - \e^TH(s)\e.
\end{equation}
Denote \(m_K = 2B^2K\left(\ln p + \sqrt{2\ln p}\right)\). From Bonferroni inequality we have
\begin{align*}
\p\left(\underset{\substack{s\colon k(s)\leq K \\ s_0 \not\subset s}}{\exists}\e^T H(s)\e \geq m_K\right) &= \p\left(\underset{1\leq j \leq K}{\exists}~\underset{\substack{s\colon k(s)=j \\ s_0 \not\subset s}}{\exists}\e^T H(s)\e \geq m_K\right)\\
& \leq \sum_{j=1}^K\binom{p}{j}\p(\e^T H(s)\e\geq m_j).
\end{align*}
From the lemma \ref{lemat2} we last sum converges to zero, so
\begin{equation*}
\underset{\substack{s\colon k(s)\leq K \\ s_0 \not\subset s}}{\forall}~\e^TH(s)\e \leq O_p(K \ln{p}),
\end{equation*}
where \(O_p(K \ln{p})\) does not depend on \(s\).

Now we prove that
\begin{equation*}
\underset{\substack{s\colon k(s)\leq K \\ s_0 \not\subset s}}{\forall}~|\mu^T[I-H(s)]\e| = \sqrt{\Delta(s) O_p(K\ln p)}.
\end{equation*}
Note that
\begin{equation}
\label{pojedyn}
\mu^T[I-H(s)]\e = \sqrt{\Delta(s)}\e'(s),
\end{equation}
where \(\e'(s)\) is a single random variable with \(B\)-subgaussian distribution. To show that, let us denote by \(\alpha\) a vector \(\mu^T[I-H(s)]\). Because \(\alpha\e\) is a sum of independent random variables with \(\alpha_iB\)-subgaussian distribution, then using the fact \ref{fakt3}, we can write that \(\alpha\e\) is \(\sqrt{\sum_{i=1}^n \alpha_i^2}B\)-subgaussian. We have
\begin{align*}
\sqrt{\sum_{i=1}^n \alpha_i^2} &= \sqrt{\alpha\alpha^T} = \sqrt{\mu^T[I-H(s)] (\mu^T[I-H(s)])^T}\\
& = \sqrt{\mu^T[I-H(s)]\mu} = \sqrt{\Delta(s)},
\end{align*}
what justifies (\ref{pojedyn}). Using the fact \ref{fakt4}, there is a constant \(c>0\) that for every \(m>0\) and fixed \(s\) we have \(\p(|\e'(s)|\geq m) \leq 2e^{-cm^2}\). Denote \(m=\sqrt{\frac{1}{c}(K+1)\ln p}\). There are \(\sum_{j=1}^K \binom{p}{j}\) sets \(s\) with the size of at most \(K\), so let us estimate
\begin{equation}
\sum_{j=1}^K \binom{p}{j}\p(|\e'(s)|\geq m) \leq 2\sum_{j=1}^K p^j \exp(-(K+1)\ln p) \leq \frac{2K}{p}.
\end{equation}
Because \(2K/p\) converges to zero, we can write
\begin{equation*}
\underset{\substack{s\colon k(s)\leq K \\ s_0 \not\subset s}}{\forall}~|\e'(s)| \leq O_p\left(\sqrt{K\ln p}\right).
\end{equation*}
Hence, we get
\begin{equation*}
\underset{\substack{s\colon k(s)\leq K \\ s_0 \not\subset s}}{\forall}~|\mu^T[I-H(s)]\e| = \sqrt{\Delta(s)O_p(K\ln p)}.
\end{equation*}

Therefore, components of the difference \(\RSS(s_0)-\RSS(s)\) are estimated in the same way like in (\ref{RSSroznica}), so we can write that \(\RSS(s)-\RSS(s_0) = \Delta(s)(1+o_p(1))\). Because
\begin{align*}
\mBIC_{\gamma}(s)-\mBIC_{\gamma}(s_0) &=  n\ln\left(1+\frac{\RSS(s)-\RSS(s_0)}{\RSS(s_0)}\right) - (2\gamma+1)p_0\ln p\\
&\geq n\ln\left(1+\frac{\Delta(s)}{n}(1+o_p(1))\right) - (2\gamma+1)p_0\ln p,
\end{align*}
then using (\ref{norm_nier}), we get that with the probability going to one \(\mBIC_{\gamma}(s)-\mBIC_{\gamma}(s_0)\) is larger than zero for every \(s\) with \(k(s)\leq K\).

\begin{center} * * * \end{center}

Now, consider the case when \(s_0\subset s\). We have \([I-H(s)]X(s_0)=0\), so
\begin{align*}
\RSS(s) &=  y^T[I-H(s)]y = \e^T[I-H(s)]\e,\\
\RSS(s_0)-\RSS(s) & = \e^T[I-H(s_0)]\e - \e^T[I-H(s)]\e = \e^T[H(s)-H(s_0)]\e
\end{align*}
and
\begin{align*}
n\ln\frac{\RSS(s_0)}{\RSS(s)} &=\ln\left(1+\frac{\RSS(s_0)-\RSS(s)}{\RSS(s_0)-[\RSS(s_0)-\RSS(s)]}\right)\\
& \leq \frac{\e^TA_j(s)\e}{\RSS(s_0)-\e^TA_j(s)\e},
\end{align*}
where \(A_j(s)\) is a symmetric idempotent matrix with rank \(j=k(s)-p_0\).

Using the lemma \ref{eAe}, we can write  that \(\RSS(s_0)= n + O_p(\sqrt n)\). Denote \(m_j = 2B^2j\left(\ln p + \sqrt{2\ln p}\right)\). From the Bonferroni inequality and the lemma \ref{lemat2} we get
\begin{align}
\label{mj2}
& \p\left(\underset{1\leq j \leq K-p_0}{\exists}~\underset{\substack{s\colon k(s)=j+p_0 \\ s_0 \subset s}}{\exists}\,\e^TA_j(s)\e \geq m_j\right)\nonumber\\
&~~~~~~~~~~~~~~~~~~~~~~~~~~~~ \leq \sum_{j=1}^K \binom{p}{j} \p(\e^TA_j(s)\e\geq m_j) \xrightarrow{n\rightarrow \infty} 0.
\end{align}

Let \(\gamma>B^2=b^2/\sigma^2\). Analogically to (\ref{subtgamma}), we can write
\begin{align}
\label{subt3}
& \p\left(\underset{\substack{s\colon k(s)\leq K \\ s_0 \subset s,\,s\neq s_0}}{\forall}~\mBIC_{\gamma}(s)-\mBIC_{\gamma}(s_0) > 0\right)\nonumber\\
& ~~~~ = 1-\p\left(\underset{\substack{s\colon k(s)\leq K \\ s_0 \subset s,\,s\neq s_0}}{\exists}~\e^T A_j(s)\e \geq (1+o_p(1))(j\ln n + 2\gamma j\ln p)\right),
\end{align}
where \(o_p(1)\) does not depend on \(s\). With the probability going to one,
\begin{equation*}
(1+o_p(1))(j\ln n + 2\gamma j\ln p)> m_j
\end{equation*}
if \(p_0\ln p=o(n)\), so from (\ref{mj2}) and (\ref{subt3}) we get the thesis.

Let \(\gamma= b^2/\sigma^2\). Analogically to (\ref{subt}), we have
\begin{align*}
& \p\left(\underset{\substack{s\colon k(s)\leq K \\ s_0 \subset s,\,s\neq s_0}}{\forall}~\mBIC_{\gamma}(s)-\mBIC_{\gamma}(s_0) > 0\right)\nonumber\\
& ~~~~ = 1-\p\left(\underset{\substack{s\colon k(s)\leq K \\ s_0 \subset s,\,s\neq s_0}}{\exists}~\e^T A_j(s)\e \geq (1+O_p(n^{\alpha-1}))(j\ln n + 2\gamma j\ln p)\right),
\end{align*}
where \(O_p(n^{\alpha-1})\) does not depend on \(s\).  With the probability going to one,
\begin{equation*}
(1+O_p(n^{\alpha-1}))(j\ln n + 2\gamma j\ln p) > m_j
\end{equation*}
if \(\sqrt{\ln p}=o(\ln n)\). 

Finally, consider \(1\leq\gamma<b^2/\sigma^2\). Then
\begin{equation*}
(1+O_p(n^{\alpha-1}))(j\ln n + 2j\ln p) > m_j
\end{equation*}
if \(p=O(n^{\delta})\) and \(2\gamma\delta+1>2\delta B^2\), so \(\delta < \frac{1}{2(b^2/\sigma^2-\gamma)}\).
\end{proof}

\begin{center} * * * \end{center}

\begin{theorem}
\label{twsub2}
Assume the model (\ref{reg}), where \(\e_i\) are \textit{b}-subgaussian with the variance \(\sigma^2\), \(p\geq n\) and let the condition (\ref{Ident}) hold. If \(p_0\ln p=o(n)\), than for \(\gamma> b^2/\sigma^2\) \(\mBIK_{\gamma}\) is consistent, that is
\begin{equation*}
\p\left(\underset{\substack{s\colon k(s)\leq K \\ s \neq s_0}}{\forall} \mBIK_{\gamma}(s)>\mBIK_{\gamma}(s_0)\right)\xrightarrow{n \rightarrow \infty} 1,
\end{equation*}
where \(K=kp_0\) for any fixed \(k>1\). If \(\gamma= b^2/\sigma^2\), \(\mBIC_{\gamma}\) is consistent if additionally \(\sqrt{\ln p}=o(\ln n)\) and \(p_0^2=O(n^{\alpha})\) for \(\alpha<1\). If \(1\leq\gamma<b^2/\sigma^2\), besides \(p_0^2=O(n^{\alpha})\) it has to be \(p=O(n^{\delta})\), where \(\delta < \frac{1-\alpha}{2(b^2/\sigma^2-1)}\).
\end{theorem}

\begin{proof}
When \(s_0 \not\subset s\), analogically to (\ref{subg3}), we can write
\begin{align*}
& \mBIK_{\gamma}(s)-\mBIK_{\gamma}(s_0)\\
& = \mBIC_{\gamma}(s) - \mBIC_{\gamma}(s_0) + 2\ln (p_0!) - 2\ln (k(s)!)\\
& \geq Cp_0\ln p \ln\left(1 + \frac{Cp_0\ln p}{n}(1+o_p(1))\right)^{\frac{n}{Cp_0\ln p}} -(2k+2\gamma+1)p_0\ln p, 
\end{align*}
so with the probability going to one \(\mBIK_{\gamma}(s)-\mBIK_{\gamma}(s_0)\) is larger than zero if \(n\) is large enough.

Now, let \(s_0\subset s\) and \(\gamma> B^2=b^2/\sigma^2\). Analogically to (\ref{subt4}), we can write
\begin{align*}
& \p\left(\underset{\substack{s\colon k(s)\leq K \\ s_0 \subset s,\,s\neq s_0}}{\forall}~\mBIK_{\gamma}(s)-\mBIK_{\gamma}(s_0) > 0\right)\\
& = 1-\p\left(\underset{\substack{s\colon k(s)\leq K \\ s_0 \subset s,\,s\neq s_0}}{\exists}~\e^T A_j(s)\e \geq (1+o_p(1))(j\ln n + 2\gamma j\ln p - 2j\ln(kp_0))\right),
\end{align*}
where \(o_p(1)\) does not depend on \(s\). We have, with the probability going to one,
\begin{equation*}
(1+o_p(1))(j\ln n + 2\gamma j\ln p - 2j\ln(kp_0))> m_j
\end{equation*}
if \(p_0\ln p=o(n)\), therefore we get the thesis.

Consider \(\gamma=b^2/\sigma^2\). Analogically to (\ref{subt2}), we get
\begin{align*}
& \p\left(\underset{\substack{s\colon k(s)\leq K \\ s_0 \not\subset s}}{\forall}~\mBIK(s)-\mBIK(s_0) > 0\right)\\
& = 1-\p\left(\underset{\substack{s\colon k(s)\leq K \\ s_0 \not\subset s}}{\exists}~\e^T A_j(s)\e \geq (1\!+\!O_p(n^{\alpha-1}))(j\ln n\!+\!2\gamma j\ln p\!-\!2j\ln(kp_0))\right),
\end{align*}
where \(O_p(n^{\alpha-1})\) does not depend on \(s\). With the probability going to one,
\begin{equation*}
(1+O_p(n^{\alpha-1}))(j\ln n + 2\gamma j\ln p -2j\ln(kp_0)) > m_j
\end{equation*}
if \(\sqrt{\ln p}=o(\ln n)\) and \(p_0^2=O(n^{\alpha})\) for \(\alpha<1\). 

Finally, consider \(1\leq\gamma<b^2/\sigma^2\). Then
\begin{equation*}
(1+O_p(n^{\alpha-1}))(j\ln n + 2j\ln p-2j\ln(kp_0)) > m_j
\end{equation*}
if \(p=O(n^{\delta})\) and \(2\gamma\delta+1>2\delta B^2 + \alpha\), so \(\delta < \frac{1-\alpha}{2(b^2/\sigma^2-\gamma)}\).
\end{proof}

\section{Simulations}
We will illustrate the theorems with simulations. Columns of the design matrix \(\textbf{X}_{n\times p}\) are generated independently from the standard normal distribution, the trait \(\textbf{y}\) according to the formula
\begin{equation*}
y_i = \sum_{j=1}^p \beta_j x_{ij} + \e_i, ~~i \in \{1,\ldots,n\}.
\end{equation*}
The number of observations, \(n\), changes in the range from 100 do 1000 every 100 and for these \(n\) we have \(p=n^2/10\) variables (so \(p\) increase from \(10^3\) to \(10^5\)), of which  \(p_0=2[n^{1/3}]\) (where \([\cdot]\) is the integer part; so \(9\leq p_0 \leq 20\)) coefficients \(\beta_j\) are equal to \(1.3n^{-1/4}\) (which gives \(0.41 \geq\beta_j\geq 0.23\)), the rest is equal to zero. The errors \(\e_i\) are generated from the standard normal and the Rademacher distribution, that is \(\p(\e_i=-1)=\p(\e_i=1)=1/2\). These parameters meet assumptions of the theorem \ref{twsub2}, because for \(\e_i\) from the Rademacher distribution we have
\begin{equation*}
\E\left(e^{t\e_i}\right)\leq \frac{1}{2}(e^t + e^{-t}) =\cosh(t)\leq e^{t^2/2},
\end{equation*}
therefore \(\e_i\) is 1-subgaussian. Besides, \(\Var(\e_i)=1\), so we can choose \(\gamma=1\). We have
\begin{align}
\label{zal}
& p_0^2 = 2[n^{2/3}] = O(n^{2/3}),\nonumber\\
& \sqrt{\ln p} = \sqrt{\ln(n^2/10)} = o(\ln n),\nonumber\\
& \sqrt{\frac{n}{p_0\ln p}}\min_{j\in s_0}|\beta_j| = 1.3\sqrt{\frac{n^{1/6}}{2\ln(n^2/10)}} \xrightarrow{n \rightarrow \infty}\infty.
\end{align} 
Additionally, according to \cite{vempala}, because elements of \(\textbf{X}\) are generated independently from the standard normal distribution, for any \(\textbf{z}\in \textbf{R}^p\) we have
\begin{equation*}
\p\left( (1-\delta)||\textbf{z}||^2 \leq ||\textbf{Xz}/\sqrt{n}||^2 \leq (1+\delta)||\textbf{z}||^2\right) \geq 1-2\exp(-(\delta^2-\delta^3)n/4).
\end{equation*}
The number of all sets \(s\) of size \(k(s)\leq kp_0\) is less than
\begin{equation*}
kp_0\binom{p}{kp_0} \leq 2kn^{1/3} \left(\frac{n^2}{10}\right)^{2n^{1/3}} = o(e^n),
\end{equation*}
then the condition (\ref{riesz}) holds with the probability going to one. As shown earlier, this condition together with (\ref{zal}) implies the consistency condition (\ref{Ident}). Because assumptions of the theorem \ref{twsub2} are the strongest, assumptions of the other theorems are also met.

Additionally, the errors from the Pareto distribution \(Pa(0,2.84)\) (that is with the density function \(f(x)=\frac{2{,}84}{x^{3{,}84}}\) for \(x\geq 1\) and zero otherwise) are generated, decreased by \(2.84/(2.84-1)\) (to get the mean equal to zero). Although the parameters are chosen so that the standard deviation is equal to one, it does not guarantee that the sample standard deviation is one, especially for a small \(n\). In fact, for \(n=100\) we get the sample standard deviation less than one in almost 80\% cases. For this reason, the error vector is divided by the sample standard deviation.

If we want to choose the best model for such a large number of explanatory variables (let us remind that for \(n = 1000\) there are \(10^5\) predictors), it is not possible to calculate criteria for each model (already for \(p = 100\) there are \(2^{100} \approx 10^{30}\) different models). What is needed is a procedure that will allow us to choose the best possible model in a reasonable time. We use the algorithm given in \cite{Prakt}, which is a modification of the stepwise procedure.

\begin{figure}
\caption{The power and FDR for errors from the normal, Rademacher and Pareto distribution after using mBIC and mBIC2. The vertical lines represent standard errors multiplied by the 0.975 quantile of the standard normal distribution.}
\label{zgodnosc}
\includegraphics{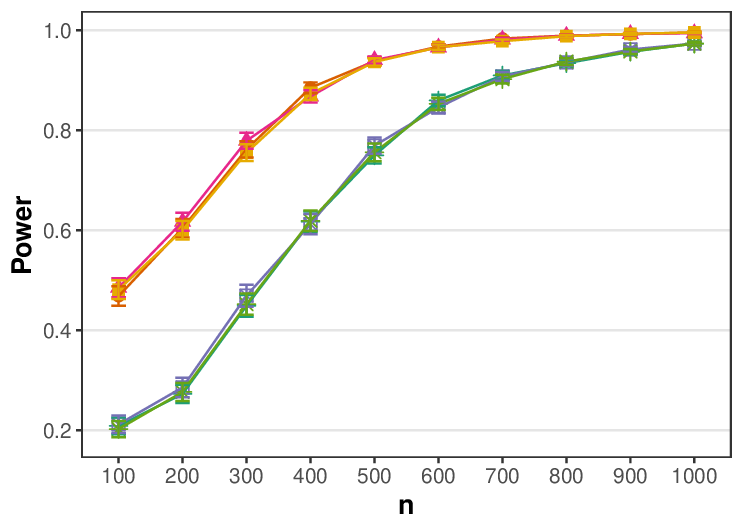}
\includegraphics{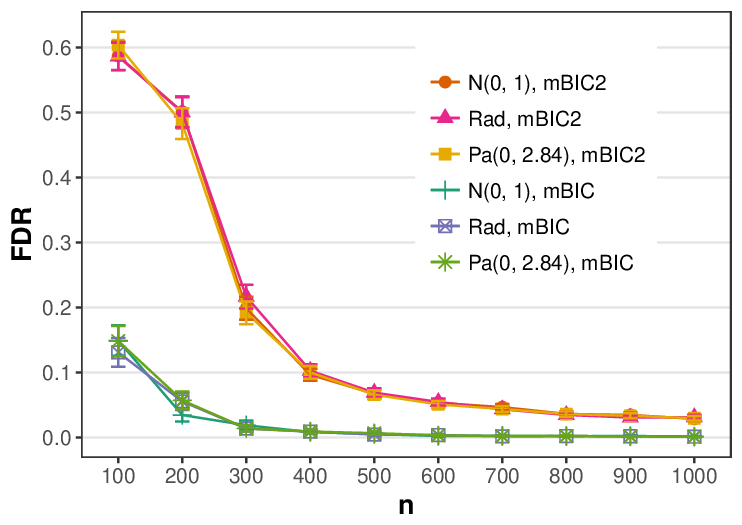}
\end{figure}

We performed \(500\) simulations. For every \(n\) the design matrix was generated once, but positions of causal columns were chosen randomly in each simulation. The false discovery rate was estimated from the formula  \(\frac{FP}{\max\{1,\,FP+TP\}}\) and the power as \(\frac{TP}{p_0}\), where \(FP\) is a number of the false positives, \(TP\) is a number of the false negatives. Results were averaged and presented in the figure \ref{zgodnosc}. Simulations were made in the R environment using the \textit{bigstep} package written by the author.

As you can see, both for the gaussian and subgaussian errors the power tends to one while FDR goes to zero, which implies the consistency. Because the penalty in mBIC2 is lower than in mBIC, the power for this criterion is higher, but FDR also. We observe similar behaviour for the Pareto distribution (which is not subgaussian), and this suggests that the theorems given in the previous section can be strengthened. For \(n=200\) FDR is greater than 0.1 (although it has been previously stated that mBIC2 controls FDR), however, we must remember that the design matrix is ​​non-orthogonal (the maximum absolute correlation between columns for the simulated data is 0.36).

\subsection{Backcross design}
\begin{figure}
\caption{The power and FDR for data simulated according to the backcross design, after using mBIC and mBIC2. The vertical lines represent standard errors multiplied by the 0.975 quantile of the standard normal distribution.}
\label{zgodnosc2}
\begin{minipage}[t]{0.47\textwidth}
\includegraphics{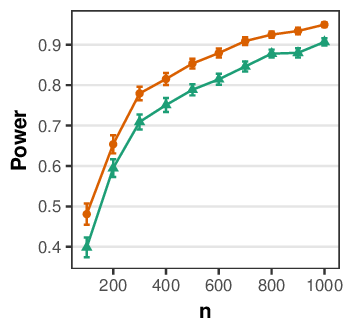}
\end{minipage}
\hspace{0.2cm}
\begin{minipage}[t]{0.47\textwidth}
\includegraphics{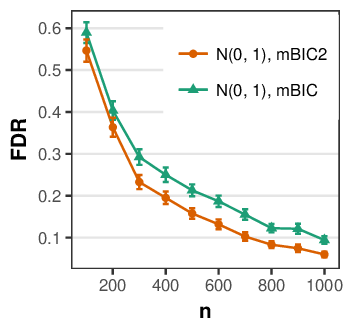}
\end{minipage}
\end{figure}

\begin{figure}
\caption{The number of true and false discoveries for data simulated according to the backcross design, after using mBIC and mBIC2. The vertical lines represent standard errors multiplied by the 0.975 quantile of the standard normal distribution.}
\label{fp}
\begin{minipage}[t]{0.47\textwidth}
\includegraphics{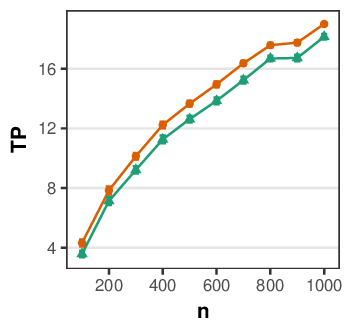}
\end{minipage}
\hspace{0.2cm}
\begin{minipage}[t]{0.47\textwidth}
\includegraphics{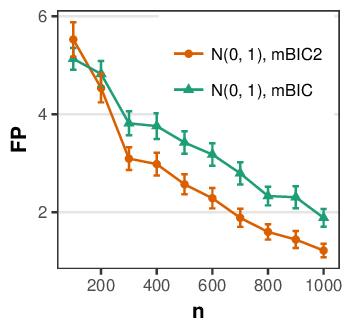}
\end{minipage}
\end{figure}

Recall that in the design matrix discussed above the elements were generated independently from the standard normal distribution. In actual genetic applications, there is often a situation in which columns of the design matrix are characterized by a high and slowly disappearing correlation, such as in the so called \textit{backcross} design. In this case, predictors can take only two values and the correlation between columns \(t\) and \(s\) is equal to \(\exp(-\alpha|t-s|)\), where \(\alpha>0\). In order to illustrate that the criteria are also consistent in this situation, the individual columns (markers) were simulated according to the backcross, for \(\alpha=0.04\). Other parameters were chosen as in previous simulations, only coefficients \(\beta_j\) were increased to \(8n^{-1/4}\), while genotypes of markers were coded as 0 and 1. The results are shown in the Figure \ref{zgodnosc2}. We see that the power goes to one, while the fraction of false discoveries goes to zero. It is interesting that this time mBIC2 has not only higher power but also lower FDR. This is not only due to the fact that mBIC2 finds more causal markers than mBIC, but also the number of false discoveries is smaller. This is illustrated by The Figure \ref{fp}, on which the numbers of true and false discoveries were given instead of the power and FDR.

This phenomenon can be explained in such a way that when a model contains more causal markers, a precision of their location increases. Because the correlation between columns is high, the criteria find markers that are very close to causal---but we count them as false discoveries. Models found by mBIC2 are larger, thanks to which a precision increases and this situation is rarer.

\subsection{Package \textit{bigstep}}
We want to add few words about the package which was used in those simulations. It is worth noting that if we have 1000 observations and \(10^5\) predictors, it is not easy to perform the stepwise procedure on such big data using home computers. What is more, it should be relatively quick because when we do simulations, we want to repeat this many times. That was a motivation to write the package \textit{bigstep}. The most important feature of this package is that it does not load a whole design matrix to the memory, which is impossible when we have a lot of predictors, but instead it keeps this matrix on the hard drive and works only on a part of it (that is the package only loads as many columns as the computer's memory allows). When we check every predictor in that part, we go to the next one. It is possible because the stepwise is a sequential procedure, we do not need to have access to every predictor at the same time. Theoretically, the package should work with any number of variables. There is only one condition: models which are fitting in next steps cannot be too large, that is they cannot exceed the memory capacity. In practise, it is very hard to do that, especially when we are interested in sparse models.

\section{Discussion}
We gave conditions when mBIC and mBIC2 are consistent, both when the error term is gaussian or subgaussian. Theorems were supported by simulations and it should be emphasized that in all of them neither the size of the true model nor the number of all variables nor the size of \(\beta_j\) were fixed (then this type of asymptotic behaviour of any reasonable criterion would be obvious). On the contrary: the size of the correct model and the number of all variables were increasing, while the magnitude of \(\beta_j\) decreased. Furthermore, unlike EBIC and a lot of other criteria, mBIC and mBIC2 are precise, that is the penalty is accurate, does not depend on constants. It is important because in practice, asymptotic behaviour is not as important as FWER or FDR. If one is interested only on building models to perform good predictions, these constants can be selected using for example cross-validation, but it does not have to be a good practice when we are interested in the inference.

However, it is worth noting that in our simulations, when looking for the best model, we consider only a small part of them. As a result, some of the observed properties may not be a result of the characteristics of the criteria, but the method used to select a subset of the analyzed models. In a future work we want to show that the whole procedure, that is mBIC/mBIC2 and the stepwise, is consistent. We believe that this is true, based on a lot of simulations we have done and results of Su \cite{su}, who showed that if the design matrix is gaussian and a signal is sparse and strong enough, the forward selection gives us the appropriate ordering of variables.

\end{document}